\documentclass[11pt,a4paper]{article}
\usepackage{amssymb, shapepar}
\usepackage{mathrsfs}
\usepackage{latexsym,bm}
\usepackage{amsthm,amsmath}
\usepackage{indentfirst}
\usepackage{color}
\usepackage{tabularx}

\theoremstyle{plain}\newtheorem{thm}{Theorem}[section]
\theoremstyle{plain}\newtheorem{defn}{Definition}[section]
\theoremstyle{plain}
\theoremstyle{plain}
\theoremstyle{plain}\newtheorem{lem}{Lemma}[section]
\theoremstyle{plain}

\newtheorem*{thmA}{Theorem A}
\newtheorem*{thmB}{Theorem B}
\newtheorem*{thmC}{Theorem C}
\newtheorem*{thmD}{Theorem D}
\newtheorem*{thmE}{Theorem E}

\numberwithin{equation}{section}

\def\rr{\mathbb{R}}
\def\cc{\mathbb{C}}
\def\rn{\mathbb{R}^{n}}

\def\zz{\mathbb{Z}}
\def\nn{\mathbb{N}}
\def\cdd{\mathcal {D}}
\def\cd'{\mathcal{D}'}
\def\cm{\mathcal {M}}

\def\az{\alpha}
\def\bz{\beta}

\def\dz{\delta}

\def\ez{\epsilon}
\def\fz{\infty}
\def\gz{\gamma}
\def\lz{\lambda}

\def\wz{\omega}
\def\pz{\partial}
\def\rz{\rho}
\def\tz{\theta}

\def\vz{\varphi}

\def\ls{\lesssim}

\def\v{$V$}
\def\B{\mathcal{B}}

\def\L{\mathcal{L}}

\def\T{\mathcal{T}}

\def\wha{$H_{\L}^{1}(\omega)$}

\def\whaf{\|f\|_{H_{\L}^{1}(\omega)}}

\def\l{\left}
\def\r{\right}

\def\supp{{\rm supp}}
\def\loc{{\rm loc}}
\def\wt{\widetilde}

\def\dsum{\displaystyle\sum}

\def\dint{\displaystyle\int}
\def\dfrac{\displaystyle\frac}
\def\dsup{\displaystyle\sup}

\def\dinf{\displaystyle\inf}

\def\fin{{\mathop\mathrm{fin}}}

\newcommand{\bea}{\begin{align}}
\newcommand{\eea}{\end{align}}
\newcommand{\beq}{\begin{equation}}
\newcommand{\eeq}{\end{equation}}
\newcommand{\beqs}{\begin{eqnarray*}}
\newcommand{\eeqs}{\end{eqnarray*}}
\newcommand{\beqn}{\begin{eqnarray}}
\newcommand{\eeqn}{\end{eqnarray}}
\newcommand{\beqa}{\begin{array}}
\newcommand{\eeqa}{\end{array}}
\newcommand{\noi}{\noindent} 

\begin{document}

\title{\bf  Riesz transform characterization of
weighted  Hardy spaces
associated to Schr\"{o}dinger operators
}

\author{ Hua Zhu }

\date{}
\maketitle

\noi {\bf Abstract.}\quad
In this paper, we characterize the weighted local Hardy spaces $h^p_\rho(\wz)$
related to the critical radius function $\rho$ and weights $\wz\in A_{1}^{\rz,\,\fz}(\rn)$
by localized Riesz transforms $\widehat{R}_j$,
in addition,  we give a characterization of  weighted
 Hardy spaces $H^{1}_{\cal L}(\wz)$ via Riesz tranforms
  associated to Schr\"{o}dinger operator $\cal L$,
where $\L=-\Delta+V$ is a  Schr\"{o}dinger operator on $\rn$ ($n\ge 3$) and
$V$ is a nonnegative function satisfying the reverse H\"older inequality.

\bigskip

\section{Introduction} \label{s1}

\noi Let $\L=-\Delta+V$ be a  Schr\"{o}dinger operator on $\rn,\ n\geq 3$,
where $V \not\equiv 0$ is a fixed non-negative potential.
We assume that $V$ belongs to the
reverse H\"{o}lder class $RH_s(\rn)$ for some $s \geq n/2$;
that is, there
exists $C = C(s,V) > 0$ such that
\beq
\l(\dfrac{1}{|B|}\dint_B V(x)^{s}\,dx \r)^{\frac{1}{s}}\leq C \l(\dfrac{1}{|B|}\int_B V(x)\,dx\r),
\eeq
for every ball $B\subset\mathbb{R}^{n}$.
Trivially, $RH_q(\rn)\subset RH_p(\rn)$ provided $1<p \leq q < \infty$.
It is well known that, if $V\in RH_q(\rn)$ for some $q>1$,
then there exists $\varepsilon>0$, which
depends only on $d$ and the constant $C$ in (1.1),
such that $V\in RH_{q+\varepsilon}(\rn)$ (see \cite{Ge}).
Moreover, the measure $V(x)\,dx$ satisfies the doubling condition:
\beqs
\dint_{B(y,2r)}V(x)\,dx\leq C\dint_{B(y,r)}V(x)\,dx.
\eeqs

With regard to the Schr\"odinger operator $\L$, we
know that the operators derived from $\L$ behave "locally" quite similar to
those corresponding to the Laplacian (see \cite{DG,Sh}).
The notion of locality is given by the critical radius function
\beq \label{2.1}
\rho(x)=\dfrac 1{m_V(x)}=\dsup_{r>0}\l\{r:\ \dfrac
{1}{r^{n-2}}\dint_{B(x,r)}V(y)dy\le 1\r\}.
\eeq

Throughout the paper we assume that $V\not\equiv0$,
so that $0<\rz(x)<\fz$ (see \cite{Sh}).
In particular, $m_V(x)=1$ with $V=1$
and $m_V(x)\sim (1+|x|)$ with $V=|x|^2$.

The study of schr\"odinger operator $\L=-\triangle+V$ recently
attracted much attention; see \cite{BHS1,BHS2,DZ1,DZ2,DZ3,LTZ,Sh,Ta3,Ta4,YZ,YZ1,Zh,ZL}.
In particularly,  J. Dziuba\'{n}ski and J. Zienkiewicz \cite{DZ1,DZ2} studied
Hardy space $H^{1}_{\cal L} $ associated to Schr\"{o}dinger operators $\cal L$ with
potential satisfying reverse H\"{o}lder inequality.
In \cite{LTZ}, we introduce the weighted Hardy space $H_{\L}^1(\wz)$
with $\wz\in A_1(\rn)$ (Muckenhoupt's weights, see \cite{Du,GR,Ga,St,ST}), and in \cite{ZL}, we obtain some weighted estimates
for bilinear operators on $H_{\L}^1(\wz)$.

Recently, Bongioanni, etc. \cite{BHS1} introduced new classes of weights,
 related to Schr\"odinger operators ${\cal L}$, that is,
 $A_p^{\rz,\fz}(\rn)$ weight
 which are in general larger than Muckenhoupt's
 (see Section 2 for notions of $A_p^{\rz,\fz}(\rn)$ weight).
In \cite{Ta2,Ta3,Ta4}, Tang gave some  results about weighted norm inequalities
with weights $\wz\in A_p^{\rz,\fz}(\rn)$.
In \cite{ZT}, we have
obtained a atomic characterization for
weighted Hardy space $H^{1}_{\cal L}(\wz)$  with $\wz\in A_1^{\rz,\fz}(\rn)$,
by introducing  the  weighted local
Hardy spaces $h^p_\rho(\wz)$ with  $A_q^{\rz,\fz}(\rn)$ weights
and establishing
the  atomic characterization of the  weighted local
Hardy spaces $h^p_\rho(\wz)$ with  $\wz\in A_q^{\rz,\fz}(\rn)$ weights.

The purpose of this paper is to give a characterization of  weighted
Hardy spaces $H^{1}_{\cal L}(\wz)$ with $\wz\in A_1^{\rz,\fz}(\rn)$
via Riesz transforms associated to Schr\"{o}dinger operator $\cal L$.

For $j=1,2,...,n$, let us define the Riesz transforms $R_j$
associated to Schr\"{o}dinger operator $\cal L$
by
$ R_j =\frac{\pz}{\pz x_j}\L^{-1/2}.$
In addition, as in \cite{YYZ}, for all $j\in \{1,2,\cdots,n\}$, and $x\in\rn$,
we define
localized Riesz transforms as
\beq
\widehat{R}_j(f)(x)\equiv \text{p.v.} \,c_n\int_{\rn}\frac{x_j-y_j}{|x-y|^{n+1}}
\eta \l(\frac{|x-y|}{\rz(x)}\r)f(y)\,dy,
\eeq
where and in what follows, $c_n\equiv \Gamma((n+1)/2)/[\pi^{(n+1)/2}]$,
$\eta\in C^1(\rn)$ supported in $(-1,1)$ and $\eta(t)=1$ if $|t|\le 1/2$.
Then, the main result of this paper is the following theorem.
(See Section 2 for definition of $h^1_{\rho}(\wz)$ and $H^1_{\L}(\wz)$.)
\begin{thm}\label{t1.1}
Let  $\wz\in A_1^{\rz,\,\fz}(\rn)$,
then we have
$$\whaf \sim\|f\|_{h^1_{\rho}(\wz)}
\sim
\|f\|_{L^1_{\wz}(\rn)}
+\sum_{j=1}^n\|\widehat{R}_j(f)\|_{L^1_{\wz}(\rn)}\sim
\|f\|_{L^1_{\wz}(\rn)}
+\sum_{j=1}^n\|{R}_j(f)\|_{L^1_{\wz}(\rn)}.$$
\end{thm}

The paper is organized as follows. In Section 2,
we review some notions and notations
concerning the weight classes $A_p^{\rz,\tz}(\rn)$
introduced in \cite{BHS1,Ta3,Ta4},
and we also review some results about $h^1_{\rho}(\wz)$ and $H^1_{\L}(\wz)$.
In Section 3,
we give the proof of Theorem 1.1.

Throughout this paper, we let $C$ denote  constants that are
independent of the main parameters involved but whose value may
differ from line to line. By $A\sim B$, we mean that there exists a
constant $C>1$ such that $1/C\le A/B\le C$.
The symbol $A\ls B$ means that $A\le CB$.
The  symbol $[s]$ for $s\in\rr$ denotes the maximal integer not more than $s$.
We also set $\nn\equiv\{1,\,2,\, \cdots\}$ and $\zz_+\equiv\nn\cup\{0\}$.
The multi-index notation is usual:
for $\az=(\az_1,\cdots,\az_n)$ and
$\pz^\az=(\pz/\pz_{x_1})^{\az_1}\cdots(\pz/\pz_{x_n})^{\az_n}$.

\section{Preliminaries}\label{s2}

\noi In this section, we first review some notions and notations
concerning the weight classes $A_p^{\rz,\tz}(\rn)$
introduced in \cite{BHS1,Ta3,Ta4}.
Given $B=B(x,r)$ and $\lz>0$, we will write $\lz B$ for
the $\lz$-dilate ball, which is the ball with
the same center $x$ and with radius $\lz r$. Similarly,
$Q(x,r)$ denotes the cube centered at $x$ with
side length $r$ (here and below only cubes with sides
parallel to the axes are considered),
and $\lz Q(x,r)=Q(x,\lz r)$.
Especially, we will denote $2 B$ by $B^*$, and $2Q$ by $Q^*$.
The following lemma is a  basic property of the critical radius function $\rz(x)$.

\begin{lem} \label{l2.1} {\bf (see \cite{Sh})} \
 There exist $C_0\ge 1$ and $k_0 \geq 1$ so that for all $x,y\in \rn$
\beq \label{2.2}
C_0^{-1}\rz(x)\l(1+\frac{|x-y|}{\rz(x)}\r)^{-k_0}\leq \rz(y)\leq C_0\rz(x)\l(1+\frac{|x-y|}{\rz(x)}\r)^{\frac{k_0}{k_0+1}}.
\eeq
 In particular, $\rz(x)\sim\rz(y)$ when $y\in B(x,r)$
 and $r\leq C\rz(x)$, where $C$ is a positive constant.
\end{lem}

A ball of the form $B(x,\rz(x))$ is called critical, and in what follows we will
call critical radius function to any positive continuous function $\rz$ that
satisfies \eqref{2.1}, not necessarily coming from a potential $V$.
Clearly, if $\rz$ is such a function, so it is $\bz \rz$ for any $\bz>0$.

Now we recall the covering of $\rn$
and the partition of unity related to $\rz$ from \cite{DZ1}.
For $m\in \zz$ we defined the sets $\B_m$ by
$\B_m = \{x:2^{-(m+1)/2}< \rz(x)\le 2^{-m/2}\},$
where critical radii $\rz(x)$ have been defined in (2.1).
Since $0<\rz(x)<\fz$, we have $\rn=\cup_{m\in\zz} \B_m$.

\begin{lem}{\bf (see \cite{DZ1})} \label{l2.2}
There is a positive constant $C$ and a collection of balls
$B_{(m,k)}=B{(x_{(m,k)},2^{2-m/2})}$,
$m\in \zz,\ k=1,2,...,$ such that $x_{(m,k)}\in \B_m,\ \B_m\subset \bigcup_k B(x_{(m,k)},2^{-m/2})$,
  and for every $(m,k)$ and $R\geq 2$
  \beqs \# \l\{(m',k'): B(x_{(m,k)},R2^{-m/2})\cap B(x_{(m',k')},R2^{-m'/2})\neq \emptyset\r\}\leq R^C,\eeqs
where and in what follows, for any set $E$, $\# E$ denotes its cardinality.
\end{lem}

\begin{lem}{\bf (see \cite{DZ1})} \label{l2.3}
There are non-negative functions $\psi_{(m,k)}$ such that
\begin{enumerate}
\item[\rm(i)] for all $m\in \zz$ and $k\in \nn$,
$\psi_{(m,k)}\in C_c^\fz(B(x_{(m,k)},2^{1-m/2}))$;
\item[\rm(ii)] $\sum_{(m,k)}\psi_{(m,k)}(x)=1$ for all $x\in\rn$;
\item[\rm(iii)] there exists a positive constant $C$ such that
for all $m\in \zz$ and $k\in \nn$,
$$ \|\nabla \psi_{(m,k)}\|_{L^{\fz}(\rn)}\le C2^m.$$
\end{enumerate}
\end{lem}

In this paper, we write $\Psi_\tz(B)=(1+r/\rho(x_0))^\tz$, where
 $\tz \ge 0$, $x_0$ and $r$ denotes the center and radius of $B$ respectively.

A weight  always refers to a positive function which is locally
integrable.
As in \cite{BHS1}, we say that a weight $\wz$ belongs to the class
$A_p^{\rho,\tz}(\rn)$ for $1<p<\fz$, if there is a constant $C$ such that for
all balls
  $B$
$$\l(\dfrac 1{\Psi_\tz(B)|B|}\dint_B\wz(y)\,dy\r)
\l(\dfrac 1{\Psi_\tz(B)|B|}\dint_B\wz^{-\frac 1{p
-1}}(y)\,dy\r)^{p-1}\le C.$$
We also
say that a  nonnegative function $\wz$ satisfies the $A_1^{\rho,\tz}(\rn)$
condition if there exists a constant $C$ such that
$$M_{V,\tz}(\wz)(x)\le C \wz(x), \ a.e.\ x\in\rn.$$
 where
$$M_{V,\tz}f(x)\equiv\dsup_{x\in B}\dfrac 1{\Psi_\tz(B)|B|}\dint_B|f(y)|\,dy.$$
When $V=0$, we
denote $M_{0}f(x)$ by $Mf(x)$ (the standard Hardy-Littlewood
maximal function). It is easy to see that $|f(x)|\le M_{V,\tz} f(x)\le
Mf(x)$ for $a.e.\ x\in\rn$ and any $\tz\ge 0$.

Clearly, the classes $A_p^{\rho,\tz}$ are increasing with $\tz$, and
we denote $A_p^{\rz,\fz}=\bigcup_{\tz\ge 0}A^{\rz,\tz}_{p}$.
By H\"older's inequality, we see that
$A^{\rz,\tz}_{p_1} \subset A^{\rz,\tz}_{p_2}$, if $1\le p_1<p_2<\fz$,
and we also denote $A_{\fz}^{\rz,\fz}=\bigcup_{p\ge 1}A^{\rz,\fz}_{p}$.
In addition, for $1\leq p\leq \fz$,
denote by $p'$ the adjoint number of $p$, i.e. $1/p+1/p'=1$.

Since $\Psi_\tz(B)\ge 1$ with $\tz\ge 0$, then
$A_p\subset A_p^{\rho,\tz}$ for $1\le p<\fz$, where $A_p$ denotes
the classical Muckenhoupt weights; see \cite{GR} and \cite{Mu}.
Moreover, the inclusions are proper.
In fact, as the example given in \cite{Ta3},
let $\tz>0$ and $0\le\gz\le\tz$,
it is easy to check that
$\wz(x)=(1+|x|)^{-(n+\gz)}\not\in A_\fz=\bigcup_{p\ge 1}A_p$ and
$\wz(x)dx$ is not a doubling measure, but
$\wz(x)=(1+|x|)^{-(n+\gz)}\in A_1^{\rho,\tz}$
provided that  $V=1$ and
$\Psi_\tz (B(x_0,r))=(1+r)^\tz$.

In what follows, given a Lebesgue measurable set $E$ and a weight $\wz$,
$|E|$ will denote the Lebesgue measure of $E$ and $\wz(E):=\int_E \wz(x)\,dx$.
For any
$\wz\in A^{\rz,\fz}_{\fz}$, the space $L^p_{\wz}(\rn)$
with $p\in(0,\fz)$ denotes the set of all measurable functions $f$
such that
$$\|f\|_{L^p_{\wz}(\rn)}\equiv\left(\int_{\rn}|f(x)|^p
\wz(x)\,dx\r)^{1/p}<\fz,$$
and $L^{\fz}_{\wz} (\rn)\equiv L^{\fz}(\rn)$.
The  symbol $L^{1,\,\fz}_{\wz}(\rn)$ denotes the
set of all measurable functions $f$ such that
$$\|f\|_{L^{1,\,\fz}_{\wz}(\rn)}\equiv\sup_{\lz>0}\left\{\lz
\wz(\{x\in\rn:\,|f(x)|>\lz\})\r\}<\fz.$$

 We remark that  balls can be replaced by cubes  in definition of
$A_p^{\rho,\tz}$ and $M_{V,\tz}$, since
$\Psi(B)\le \Psi(2B)\le 2^\tz \Psi(B).$
In fact, for the cube $Q=Q(x_0,r)$, we can also define
$\Psi_\tz(Q)=(1+r/\rho(x_0))^\tz$.
Then  we give the weighted boundedness of $M_{V,\tz}$.

Now, we give some properties of weights class $A^{\rz,\tz}_{p}$ for $p\ge1$.

\begin{lem}{\bf (see \cite{BHS1,Ta4,ZT})} \label{l2.4}
Let $\wz\in A^{\rho,\fz}_p=\bigcup_{\tz\ge 0}A_p^{\rho,\tz}$ for $p\ge 1$.
Then
\begin{enumerate}
 \item[\rm(i)] If $ 1\le p_1<p_2<\fz$, then $A_{p_1}^{\rho,\tz}\subset
A_{p_2}^{\rho,\tz}$.

 \item[\rm(ii)] $\wz\in A_p^{\rho,\tz}$ if and only
if\ $\wz^{-\frac 1{p-1}}\in A_{p'}^{\rho,\tz}$, where $1/p+1/p'=1.$

 \item[\rm(iii)] If $\wz\in A_p^{\rho,\fz},\ 1<p<\fz$, then there exists
$\ez>0$ such that $\wz\in A_{p-\ez}^{\rho,\fz}.$

 \item[\rm(iv)]Let $f\in L_{loc}(\rz)$, $0<\dz<1$, then
$(M_{V,\tz}f)^\dz\in A_1^{\rho,\tz}$.

 \item[\rm(v)] Let $1<p<\fz$, then $\wz\in A_p^{\rho,\fz}$ if and only
if $\wz=\wz_1\wz_2^{1-p}$, where $\wz_1,\wz_2\in A_1^{\rho,\fz}$.

 \item[\rm(vi)] For $\wz\in A_p^{\rho,\tz}$,  $Q=Q(x,r)$ and $\lz>1$,  there
 exists a positive constant $C$ such that
 $$\wz(\lz Q)\le C (\Psi_{\tz}(\lz Q))^p\, \lz^{np}\, \wz(Q).$$

\item[\rm(vii)] If $p\in(1,\fz)$ and $\wz\in A^{\rz,\tz}_{p}(\rn)$,
then the local Hardy-Littlewood maximal operator $M^{\loc}$ is
bounded on $L^p_{\wz}(\rn)$.

\item[\rm(viii)] If $\wz\in A^{\rz,\tz}_{1}(\rn)$, then $M^{\loc}$ is
bounded from $L^1_{\wz}(\rn)$ to $L^{1,\,\fz}_{\wz}(\rn)$.

\end{enumerate}

\end{lem}

For any $\wz\in A^{\rz,\fz}_{\fz}(\rn)$, define the  critical
index of $\wz$ by
\beq \label{qw}
q_{\wz}\equiv\inf\left\{p\in[1,\fz):\,\wz\in A^{\rz,\fz}_p(\rn)\r\}.
\eeq
Obviously, $q_\wz\in [1,\fz)$. If $q_\wz\in(1,\fz)$, then
$\wz\not\in A_{q_\wz}^{\rz,\fz}$.

The symbols ${\cal D}(\rn)=C_0^\fz(\rn), {\cal D}'(\rn)$ is the
dual space of ${\cal D}(\rn)$, and for ${\cal D}(\rn),\,{\cal D}'(\rn)$
and $L^p_{\wz}(\rn)$,
For any $\vz\in {\cal D}(\rn)$, let
$\vz_t(x)=t^{-n}\vz\l(x/t\r)$ for $t>0$ and
$\vz_j(x)=2^{jn}\vz\l(2^jx\r)$ for $j\in \zz$.
It is easy to see that we have the following results.

Next let us introduce the space $H^1_{\L}(\wz)$.

Let $\{T_t\}_{t>0}$ be the semigroup of linear operators generated by $\L$ and
$T_t(x,y)$ be their kernels, that is,
\beq
T_t f(x)=e^{-t\L}f(x)=\int_{\rn}T_t(x,y)f(y)\,dy,
\qquad {\text{for}}\ t>0\ {\text{and}}\ f\in L^{2}(\rn).
\eeq
Since \v \ is non-negative the Feynman-Kac formula implies that
\beq
0\leq T_t(x,y)\leq \wt{T}_t(x,y)\equiv(4\pi t)^{-\frac{n}{2}}\exp\l(-\frac{|x-y|^{2}}{4t}\r).
\eeq
Obviously, by (2.4) the maximal operator
\beqs
\T^{*}f(x)=\dsup_{t>0}\l|T_tf(x)\r|
\eeqs
is of weak-type (1,1).
A weighted Hardy-type space related to $\L$ with $A_1^{\rz,\tz}(\rn)$ weights
is naturally defined by:
\beq
H_{\L}^{1}(\omega)\equiv\{f\in L^{1}_\wz(\rn): \T^{*}f(x)\in L^{1}_\wz(\rn)\},
\qquad\text{with}
\quad \whaf \equiv\|\T^{*}f\|_{L^{1}_\wz(\rn)}.
\eeq
The $H_{\L}^{1}(\omega)$ with $\wz\in A_1(\rn)$ has been studied in \cite{LTZ,ZL},
and we have following theorem.

\begin{thmA} {\bf (see \cite{LTZ})} \
If $\wz\in A_1(\rn)$ and $V\in RH_{n}(\rn)$ is a non-negative potential,
then there is a constant $C>0$
such that for all $f\in H_{\L}^1(\wz)$
$$ C^{-1}\whaf \leq \|f\|_{L^1_\wz(\rn)}+\dsum_{j=1}^n\|R_j f\|_{L^1_\wz(\rn)}
\leq C\whaf. $$
\end{thmA}

Then let us recall some results about weighted local Hardy spaces $h^p_{\rz}(\wz)$.

We first introduce some local maximal functions. For $N\in\zz_+$
and $R\in(0,\fz)$, let
$$\begin{array}{cl}
\cdd_{N,\,R}(\rn)\equiv \Bigg\{\vz\in\cdd(\rn):&\,\supp(\vz)
 \subset B(0,R), \\
&\left.\|\vz\|_{\cdd_{N}(\rn)}\equiv\dsup_{x\in\rn}
\dsup_{{\az\in\zz^n_+},\,{|\az|\le N}}
|\partial^{\az}\vz(x)|\le1\r\}.
\end{array}$$

\begin{defn}\label{d3.1}
Let $N\in\zz_+$ and $R\in(0,\fz)$. For any $f\in\cd'(\rn)$, the
local nontangential grand maximal function $\wt{\cm}_{N,\,R} (f)$
of $f$ is defined by setting, for all $x\in\rn$,
\beq \label{tgnr}
\wt{\cm}_{N,\,R} (f)(x)\equiv\sup\left\{|\vz_l \ast
f(z)|:\,|x-z|<2^{-l}<\rz(x),\,\vz\in\cdd_{N,\,R}(\rn)\r\},
\eeq
and the  local vertical grand maximal function $\cm_{N,\,R}(f)$
of $f$ is defined by setting, for all $x\in\rn$,
\beq \label{gnr}
{\cm}_{N,\,R} (f)(x)\equiv\sup\left\{|\vz_l \ast
f(x)|:\,0<2^{-l}< \rz(x),\,\vz\in\cdd_{N,\,R}(\rn)\r\},
\eeq

\end{defn}

For convenience's sake, when $R=1$, we denote $\cdd_{N,\,R}(\rn)$,
$\wt{\cm}_{N,\,R} (f)$ and $\cm_{N,\,R}(f)$ simply by
$\cdd^0_{N}(\rn)$, $\wt{\cm}^0_{N}(f)$ and $\cm^0_{N}(f)$,
respectively;
when $R=\max\{ R_1,\ R_2,\ R_3 \}>1$
(where $R_1$, $R_2$ and $R_3$  are defined as in Lemma 4.2, 4.4 and 4.8 in \cite{ZT}),
we denote $\cdd_{N,\,R}(\rn)$,
$\wt{\cm}_{N,\,R} (f)$ and $\cm_{N,\,R}(f)$ simply by
$\cdd_{N}(\rn)$, $\wt{\cm}_{N}(f)$ and $\cm_{N}(f)$, respectively.
For any $N\in\zz_+$ and $x\in\rn$, obviously,
$$\cm^{0}_N(f)(x)\le\cm_N (f)(x)\le\wt{\cm}_N (f)(x).$$

Now we introduce the weighted local Hardy space via the local
grand maximal function as follows.

\begin{defn}\label{d3.2}
Let $\wz\in A^{\rz,\fz}_{\fz}(\rn)$, $q_{\wz}$ be as in \eqref{qw}, $p\in(0,1]$
and $\wt{N}_{p,\wz}\equiv[n(\frac {q_\wz}p-1)]+2.$
For each $N\in\nn$ with $N\ge\wt{N}_{p,\,\wz}$,
the  weighted local Hardy space is defined by
$$h^{p}_{\rho,\,N}(\wz)\equiv\left\{f\in\cd'(\rn):\ \cm_N (f)\in
L^{p}_{\wz}(\rn)\r\}.$$
Moreover, let
$\|f\|_{h^{p}_{\rho,\,N}(\wz)}\equiv\|\cm_N(f)\|_{L^{p}_{\wz}(\rn)}$.

\end{defn}

Obviously, for any integers $N_1$ and $N_2$
with $N_1\ge N_2\ge \wt{N}_{p,\,\wz}$,
$$h^{p}_{\rho,\,\wt{N}_{p,\wz}}(\wz)\subset
h^{p}_{\rho,\,N_2}(\wz) \subset h^{p}_{\rho,\,N_1}(\wz),$$
and the inclusions are continuous.

Next, we introduce some local vertical, tangential and nontangential
maximal functions, and then we establish the characterizations of the
weighted local Hardy space $h^{p}_{\rho,\,N}(\wz)$ by
these local maximal functions.

\begin{defn}\label{d3.5}
Let
\beq \label{3.3}
\psi_0\in\cdd(\rn)\,\, \text{with}\,\,\int_{\rn}\psi_0 (x)\,dx\neq0.
\eeq
For every $x\in\rn$, there exists an integer $j_x\in\zz$ satisfying
$2^{-j_x}< \rz(x)\le 2^{-j_x+1}$,
and then for $j\ge j_x$, $A,\,B\in[0,\fz)$ and $y\in\rn$, let
$m_{j,\,A,\,B,\,x}(y)\equiv(1+2^j |y|)^A 2^{{B|y|}/{\rz(x)}}$.

The  local vertical maximal function $\psi_0^{+}(f)$ of $f$ associated to
$\psi_0$ is defined by setting, for all $x\in\rn$,
\beq \label{3.4}
\psi_0^{+}(f)(x)\equiv \sup_{j\ge j_x}|(\psi_0)_j \ast f(x)|,
\eeq
the local tangential Peetre-type  maximal function
$\psi^{\ast\ast}_{0,\,A,\,B}(f)$ of $f$ associated to $\psi_0$ is
defined by setting, for all $x\in\rn$,
\beq \label{3.5}
\psi^{\ast\ast}_{0,\,A,\,B}(f)(x)\equiv\sup_{j\ge j_x,\,y\in\rn}
\frac{|(\psi_0)_j \ast f(x-y)|}{m_{j,\,A,\,B,\,x}(y)},
\eeq
and the  local nontangential maximal function
$(\psi_0)^{\ast}_{\triangledown}(f)$ of $f$ associated to $\psi_0$ is
defined by setting, for all $x\in\rn$,
\beq \label{3.6}
(\psi_0)^{\ast}_{\triangledown}(f)(x)\equiv
\sup_{|x-y|<2^{-l}<\rz(x)}|(\psi_0)_l \ast f(y)|,
\eeq
where $l\in\zz$.

\end{defn}

Obviously, for any $x\in\rn$, we have
$$\psi_0^{+}(f)(x)\le(\psi_0)^{\ast}_{\triangledown}(f)(x)
\ls \psi^{\ast\ast}_{0,\,A,\,B}(f)(x).$$
We point out that the local
tangential Peetre-type maximal function
$\psi^{\ast\ast}_{0,\,A,\,B}(f)$ was introduced by Rychkov
\cite{Ry}.

We have  obtained  the
local vertical and the local nontangential maximal function
characterizations of $h^{p}_{\rz,\,N}(\wz)$ with $N\ge N_{p,\,\wz}$ as follows.
Here and in what follows,
\begin{equation}\label{2.12}
N_{p,\,\wz}\equiv\max\left\{\wt{N}_{p,\,\wz},\,N_0\r\},
\end{equation}
where $\wt{N}_{p,\,\wz}$ and $N_0$ are  respectively as in
Definition 2.2 and Theorem 3.1 of \cite{ZT}.

\begin{thmB} {\bf (see \cite{ZT})} \
Let  $\wz\in A^{\rz,\,\fz}_{\fz}(\rn)$, $\psi_0$ and $N_{p,\,\wz}$ be respectively
as in \eqref{3.3} and \eqref{2.12}. Then for any integer $N\ge N_{p,\,\wz}$,
the following are equivalent:
\begin{enumerate}
\item[\rm(i)] $f\in h^{p}_{\rz,\,N}(\wz);$
\item[\rm(ii)] $f\in\cd'(\rn)$ and $\psi^{+}_0 (f)\in L^{p}_{\wz}(\rn);$
\item[\rm(iii)] $f\in\cd'(\rn)$ and $(\psi_0)^{\ast}_{\triangledown}
(f)\in L^{p}_{\wz}(\rn);$
\item[\rm(iv)] $f\in\cd'(\rn)$ and $\wt{\cm}_N (f)\in L^{p}_{\wz}(\rn);$
\item[\rm(v)] $f\in\cd'(\rn)$ and $\wt{\cm}^0_N (f)\in L^{p}_{\wz}(\rn);$
\item[\rm(vi)] $f\in\cd'(\rn)$ and $\cm^0_N (f)\in L^{p}_{\wz}(\rn)$.
\end{enumerate}
Moreover, for all $f\in h^p_{\rz,\,N}(\wz)$
\begin{align}\label{3.30}
\|f\|_{h^{p}_{\rho,\,N}(\wz)}&\sim \left\|\psi^{+}_0(f)\r\|_{L^{p}_{\wz}(\rn)}
\sim\left\|(\psi_0)^{\ast}_{\triangledown}(f)\r\|_{L^{p}_{\wz}(\rn)}\nonumber\\
&\sim \left\|\wt{\cm}_N(f)\r\|_{L^{p}_{\wz}(\rn)}
\sim\left\|\wt{\cm}^0_N(f)\r\|_{L^{p}_{\wz}(\rn)}
\sim\left\|\cm^0_N(f)\r\|_{L^{p}_{\wz}(\rn)},
\end{align}
where the implicit constants are independent of $f$.
\end{thmB}

Next, we introduce the weighted local atoms, via which,
we give the definition of
the weighted atomic local Hardy space.

\begin{defn}\label{d3.3}
Let $\wz\in A^{\rz,\fz}_{\fz}(\rn)$, $q_{\wz}$ be as in \eqref{qw}.
A triplet $(p,q,s)_\wz$ is called to be admissible, if $p\in (0,1]$,
$q\in(q_\wz,\fz]$ and $s\in\nn$ with $s\ge [n(q_\wz/p-1)]$.
A function $a$ on $\rn$ is said to be a $(p,q,s)_\wz-atom$ if
\begin{enumerate}
\item[\rm(i)] $\supp \, a\subset Q(x,r)$ and $r\le L_{1}\rz(x)$,
\item[\rm(ii)] $\|a\|_{L^q_\wz(\rn)}\le [\wz( Q)]^{1/q-1/p}$,
\item[\rm(iii)]$\dint_{\rn} a(x)x^\az dx=0$ for all $\az\in \zz_+^n$ with
$|\az|\le s$, when $Q=Q(x,r),\ r<L_{2}\rz(x)$,
\end{enumerate}
where $L_1=4 C_0 (3\sqrt{n})^{k_0}$, $L_2=1/C_0^2(3\sqrt{n})^{k_0+1}$,
and $C_0$, $k_0$ are constant given in Lemma 2.1.
 Moreover, a function $a(x)$ on $\rn$ is called  a
$(p,\,q)_{\wz}$-single-atom with $q\in(q_{\wz},\fz]$, if
$$\|a\|_{L^q_{\wz}(\rn)}\le[\wz(\rn)]^{1/q-1/p}.$$
\end{defn}

\begin{lem}\label{l2.5} {\bf (see \cite{ZT})} \
Let $\wz\in A_q^{\rz,\,\tz}(\rn)$ and $a$ be a $(p,q,s)_\wz$-atom,
which satisfies $\supp\, a\subset Q(x_0,r)$,
then there exists a constant $C$ such that:
$$\|a\|_{L^1(\rn)}\le C |Q| \wz(Q)^{-1/p}\Psi_{\tz}(Q).$$
\end{lem}

\begin{defn}\label{d3.4}
Let $\wz\in A^{\rz,\fz}_{\fz}(\rn)$, $q_{\wz}$ be as in \eqref{qw},
and $(p,q,s)_\wz$ be admissible, The  weighted atomic local Hardy space
$h^{p,\,q,\,s}_{\rho}(\wz)$ is defined as the set of all
$f\in\cd'(\rn)$ satisfying that
$$f=\sum_{i=0}^{\fz}\lz_i a_i$$
in $\cd'(\rn)$, where $\{a_i\}_{i\in\nn}$ are
$(p,\,q,\,s)_{\wz}$-atoms with $\supp(a_i)\subset Q_i$, $a_0$ is a
$(p,\,q)_{\wz}$-single-atom, $\{\lz_i\}_{i\in\zz_+}\subset\cc$.
Moreover, the quasi-norm of $f\in h^{p,q,s}_{\rho}(\wz)$ is defined by
$$\|f\|_{h^{p,q,s}_{\rho}(\wz)}\equiv
\dinf\l\{\l[\dsum_{i=0}^\fz|\lz_i|^p\r]^{1/p}\r\},$$
where the infimum is taken over all the decompositions of $f$ as
above.
\end{defn}

It is easy to see that if  triplets $(p,q,s)_\wz$ and $(p,\bar
q,\bar s)_\wz$ are admissible and  satisfy $\bar q\le q$ and $\bar
s\le s$, then $(p,q,s)_\wz$-atoms are $(p,\bar q,\bar
s)_\wz$-atoms, which  implies that
$h^{p,q,s}_{\rho}(\wz)\subset h^{p,\bar q,\bar s}_{\rho}(\wz)$ and
the inclusion is continuous.

For $h^p_{\rho, N}(\wz)$ and $h^{p,q,s}_{\rho}(\wz)$,
we have proved following theorem, and in virtue of the theorem, we can simplicity
denote by $h^p_{\rho}(\wz)$ the
weighted local Hardy space $h^p_{\rho, N}(\wz)$  when $N\ge N_{p,\,\wz}$.

\begin{thmC} {\bf (see \cite{ZT})} \
Let $\wz\in A_\fz^{\rz,\,\fz}(\rn)$, $q_\wz$ and $N_{p,\,\wz}$
be respectively as in \eqref{qw} and \eqref{2.12}.
If $q\in (q_\wz,\fz],\ p\in (0,1]$,  and
integers $s$ and $N$ satisfy $N\ge N_{p,\wz}$ and
$N>s\ge[n(q_\wz/p-1)]$, then $h^{p,q,s}_{\rho}(\wz)=h^p_{\rho, N}(\wz)=
h^p_{\rho, N_{p,\wz}}(\wz)$ with equivalent norms.
\end{thmC}

\begin{defn}
Let  $\wz\in A^{\rz,\,\fz}_{\fz}(\rn)$ and
$(p,\,q,\,s)_{\wz}$ be admissible as in Definition \ref{d3.3}.
Then  $h^{p,\,q,\,s}_{\rho,\,\fin}(\wz)$ is defined
to be the vector space of all finite linear combinations
of $(p,\,q,\,s)_{\wz}$-atoms and
a $(p,\,q)_{\wz}$-single-atom, and the norm of $f$
in $h^{p,\,q,\,s}_{\rho,\,\fin}(\wz)$ is defined by
\beqs
\begin{aligned}
\|f\|_{h^{p,q,s}_{\rho,\fin}(\wz)}
&\equiv\dinf\Big\{\Big[\dsum_{i=0}^k|\lz_i|^p\Big]^{1/p}:
f=\dsum_{i=0}^k\lz_ia_i,\ k\in\zz_+,\ \{\lz_i\}_{i=0}^k\subset\cc,\
\{a_i\}_{i=1}^k\ {\rm are}\ \Big. \\
&\qquad\qquad (p,q,s)_\wz\ {\rm atoms},\ {\rm and}\ a_0\ {\rm is\ a}\ (p,q)_\wz\ {\rm single\ atom}\Big\}.
\end{aligned}
\eeqs
\end{defn}

 As an application of finite atomic
decompositions, we establish boundedness in $h^p_\rz(\wz)$ of
quasi- Banach-valued sublinear operators.

As in \cite{BLYZ}, a  quasi-Banach space space ${\cal B}$
is a vector space endowed with a quasi-norm $\|\cdot\|_{{\cal B}}$
which is nonnegative, non-degenerate (i.e., $\|f\|_{\cal B}=0$ if
and only if $f=0$), homogeneous, and obeys the quasi-triangle
inequality, i.e., there exists a positive constant $K$ no less
than $1$ such that for all $f,g\in{\cal B}$, $\|f+g\|_{{\cal
B}}\le K(\|f\|_{\cal B}+\|g\|_{\cal B})$.

Let $\bz\in (0,1]$. A quasi-Banach space ${\cal B}_\bz$ with the
quasi-norm $\|\cdot\|_{{\cal B}_\bz}$ is called a
$\bz$-quasi-Banach space if $\|f+g\|^\bz_{{\cal B}_\bz}\le
\|f\|^\bz_{{\cal B}_\bz}+\|g\|^\bz_{{\cal B}_\bz}$ for all $f,g\in
{\cal B}_\bz$.

Notice that any Banach space is a $1$-quasi-Banach space, and the
quasi-Banach space $l^\bz,\ L^\bz_\wz(\rn)$ and $h^\bz_\wz(\rn)$
with $\bz\in (0,1)$ are typical $\bz$-quasi-Banach spaces.

For any given $\bz$-quasi-Banach space ${\cal B}_\bz$ with $\bz\in
(0,1]$ and  a linear space ${\cal Y}$, an operator $T$ from
 ${\cal Y}$ to ${\cal B}_\bz$ is said to be ${\cal B}_\bz$-sublinear
 if for any $f, g\in {\cal B}_\bz$  and $\lz,\ \nu\in\cc$,
 $$\|T(\lz f+\nu g)\|_{{\cal B}_\bz}\le\l(|\lz|^\bz\|T(f)\|_{{\cal
 B}_\bz}^\bz +|\nu|^\bz\|T(g)\|_{{\cal B}_\bz}^\bz\r)^{1/\bz}$$
and
$\|T(f)-T(g)\|_{{\cal B}_\bz}\le \|T(f-g)\|_{{\cal B}_\bz}$.

We remark that if $T$ is linear, then it is ${\cal B}_\bz$-sublinear.
Moreover, if ${\cal B}_\bz$ is a space of functions, and
$T$ is nonnegative and sublinear in the classical sense, then $T$
is also ${\cal B}_\bz$-sublinear.

\begin{thmD} {\bf (see \cite{ZT})}\
Let $\wz\in A_\fz^{\rz,\,\fz}(\rn), 0<p\le\bz\le 1$, and ${\cal B}_\bz$ be a
$\bz$-quasi-Banach space. Suppose $q\in (q_\wz,\fz)$ and
$T: h^{p,q,s}_{\rho,\,\fin}(\wz)\to {\cal B}_\bz$
is a ${\cal B}_\bz$-sublinear operator  such that
$$
S\equiv\sup\{\|T(a)\|_{{\cal B}_\bz}:
\ a \ {\rm is\ a\ }
(p,q,s)_\wz{\ \rm atom}
\ {\rm or}\ (p,q)_\wz\ {\rm  single\ atom}\}<\fz.
$$
Then there exists a unique bounded ${\cal B}_\bz$-sublinear
operator $\wt T$ from $h^p_\rho(\wz)$ to ${\cal B}_\bz$ which extends $T$.
\end{thmD}

Finally, in \cite{ZT}, we have
applied the  atomic characterization of the  weighted local
Hardy spaces $h_{\rho}^1(\wz)$ with  $A_1^{\rz,\tz}(\rn)$ weights  to establish atomic characterization of weighted Hardy space \wha\
associated to Schr\"{o}dinger operator  with $A_1^{\rz,\tz}(\rn)$ weights.

\begin{thmE} {\bf (see \cite{ZT})}\
Let $0\not\equiv V\in RH_{n/2}$ and $\wz\in A_1^{\rz,\,\fz}(\rn)$,
then $h^1_{\rho}(\wz)=H^1_{\L}(\wz)$ with equivalent norms,
that is
$$\|f\|_{h^1_{\rho}(\wz)}\sim\|f\|_{H^1_{\L}(\wz)}.$$
\end{thmE}

\section{Proof of Theorem 1.1}

The main purpose of this section is to give the proof of Theorem 1.1,
and we begin with some useful lemmas.

\begin{lem} {\bf (see \cite{LTZ})} \label{l3.1}
If $\wz\in A_1(\rn)$, then
there is a constant $C>0$ such that
\beqs \sum_{(m,k)}\|R_j(\psi_{(m,k)}f)-\psi_{(m,k)}R_j(f)\|_{L^1(\wz)}\le
C\|f\|_{L^1(\wz)}.\eeqs
\end{lem}

As in \cite{Ta1}, we can obtain the following lemma. Its proof is similar to
the Lemma 8.2 in \cite{Ta1}, and we omit the details here.

\begin{lem}\label{l3.2}
Let  $\widehat{R}_j$ be as above, then
\begin{enumerate}
\item[\rm(i)] $\|\widehat{R}_j(f)\|_{L^p_{\wz}(\rn)}\le C_{p,\,\wz}\|f\|_{L^p_{\wz}(\rn)}$,
    for $1<p<\fz$ and $\wz\in A_p^{\rz,\,\fz}(\rn)$,
\item[\rm(ii)] $\|\widehat{R}_j(f)\|_{L^{1,\,\fz}_{\wz}(\rn)}\le C_{\wz}\|f\|_{L^1_{\wz}(\rn)}$,
    for  $\wz\in A_1^{\rz,\,\fz}(\rn)$.

\end{enumerate}

\end{lem}

Now let us state the proof of Theorem 1.1.

{\bf {The proof of Theorem 1.1:}}

$\whaf \sim\|f\|_{h^1_{\rho}(\wz)}$ is obvious by Theorem E.
We now prove $\|f\|_{h^1_{\rho}(\wz)}
\sim
\|f\|_{L^1_{\wz}(\rn)}
+\sum_{j=1}^n\|\widehat{R}_j(f)\|_{L^1_{\wz}(\rn)}$. We first assume that $\|f\|_{L^1_{\wz}(\rn)}
+\sum_{j=1}^n\|\widehat{R}_jf\|_{L^1_{\wz}(\rn)}<\fz$, and here we
will borrow some idea from \cite{MPR}.
Let $B_{(m,k)}$ be a ball as in Lemma 2.2, $\zeta_{a_1B_{(m,k)}}$
be a $C_0^{\fz}$ nonnegative function supported in $2a_1B_{(m,k)}$
and $\zeta_{a_1B_{(m,k)}}=1$ on $a_1B_{(m,k)}$, where $a_1=3C_0+1$.
By Lemma 2.4, we can set $\bar{\wz}\in A_1(\rn)$ so that
$\bar{\wz}=\wz$ on $a_5B_{(m,k)}$,
where $a_5$ is a constant independent of  $B_{(m,k)}$
and will be given in the following proof.
Taking $\varphi\in\mathcal{D}^0_N(\rn)$, by  Theorem A,
we have
\begin{align}
\Big\|\sup_{0<2^{-l}<\rz(x)}|\varphi_l*f|\Big\|_{L^1_{\wz}(B_{(m,k)})}
&\lesssim
\Big\|\sup_{0<2^{-l}<\rz(x)}
\l|\varphi_l*(f\zeta_{a_1B_{(m,k)}})\r|\Big\|_{L^1_{\bar{\wz}}(\rn)}\nonumber\\
&\lesssim \l\|f\zeta_{a_1B_{(m,k)}}\r\|_{H^1_{\cal L}(\bar{\wz})}\nonumber\\
&\lesssim \l\|f\zeta_{a_1B_{(m,k)}}\r\|_{L^1_{\bar{\wz}}(\rn)}
+\sum_{j=1}^n\l\|R_j(f\zeta_{a_1B_{(m,k)}})\r\|_{L^1_{\bar{\wz}}(\rn)}.
\end{align}
Let us denote  the integral kernel of the operator $(\pz/\pz\, x_j)\L^{-1/2}$
by $R_j(x,y)$, and let $\tilde{R}_j(x,y)$ be the kernel of the classical Riesz
transform $\tilde{R}_jf=(\pz/\pz\, x_j)\Delta^{-1/2}f$.
Then we can define
$$\overline{R}_j(f)(x)\equiv\int_{\rn}R_j(x,y)\eta \l(\frac{|x-y|}{\rz(x)}\r)f(y)\,dy,$$
and we also have
$$\widehat{R}_j(f)(x)=
\int_{\rn}\tilde{R}_j(x,y)\eta \l(\frac{|x-y|}{\rz(x)}\r)f(y)\,dy.$$
Let $a_2=C_0(1+16a_1)+2a_1<a_5$, then we have
\beqs
\begin{aligned}
&\l\|R_j(f\zeta_{a_1B_{(m,k)}})
-\zeta_{a_1B_{(m,k)}}\widehat{R}_j(f)\r\|_{L^1_{\bar{\wz}}(\rn)}\\
&\lesssim
\l\|R_j(f\zeta_{a_1B_{(m,k)}}\chi_{a_2B_{(m,k)}})
-\zeta_{a_1B_{(m,k)}}{R}_j(f\chi_{a_2B_{(m,k)}})\r\|_{L^1_{\bar{\wz}}(\rn)}\\
&\quad +\l\|\zeta_{a_1B_{(m,k)}}R_j(f\chi_{a_2B_{(m,k)}})
-\zeta_{a_1B_{(m,k)}}\overline{R}_j(f\chi_{a_2B_{(m,k)}})\r\|_{L^1_{\bar{\wz}}(\rn)}\\
&\quad +\l\|\zeta_{a_1B_{(m,k)}}\overline{R}_j(f\chi_{a_2B_{(m,k)}})
-\zeta_{a_1B_{(m,k)}}\widehat{R}_j(f\chi_{a_2B_{(m,k)}})\r\|_{L^1_{\bar{\wz}}(\rn)}\\
&\equiv I_1+I_2+I_3.
\end{aligned}
\eeqs
For $I_1$, by the same method of the proof of Lemma 3.1 (Lemma 3.14 of \cite{LTZ}),
we get
\beq
I_1\lesssim \l\|f\r\|_{L^1(a_2B_{(m,k)},\,\bar{\wz})}
\lesssim \l\|f\r\|_{L^1_{\wz}(a_2B_{(m,k)})}.
\eeq
For $I_2$, by Lemma 3 of \cite{BHS0}, we have
\beqs
\begin{aligned}
I_2&\lesssim\l\|\zeta_{a_1B_{(m,k)}}R_j(f\chi_{a_2B_{(m,k)}})
-\zeta_{a_1B_{(m,k)}}\overline{R}_j(f\chi_{a_2B_{(m,k)}})\r\|_{L^1_{\bar{\wz}}(\rn)}\\
&=\l\|R_j(f\chi_{a_2B_{(m,k)}})
-\overline{R}_j(f\chi_{a_2B_{(m,k)}})\r\|_{L^1(2a_1B_{(m,k)},\,\bar{\wz})}\\
&\lesssim \int_{2a_1B_{(m,k)}}\int_{|x-y|>\rz(x)/2}
|R_j(x,y)|\l|1-\eta\l(\frac{|x-y|}{\rz(x)}\r)\r|
|f\chi_{a_2B_{(m,k)}}(y)|\,dy\,\bar{\wz}(x)\,dx\\
&\lesssim \int_{2a_1B_{(m,k)}}\int_{|x-y|>\rz(x)/2}
\frac{|f\chi_{a_2B_{(m,k)}}(y)|}{|x-y|^n}
\,dy\,\bar{\wz}(x)\,dx\\
&\lesssim \int_{2a_1B_{(m,k)}}\int_{|x-y|>\rz(x)/2}
\frac{|f\chi_{a_2B_{(m,k)}}(y)|}{(\rz(x))^n}
\,dy\,\bar{\wz}(x)\,dx.\\
\end{aligned}
\eeqs
Since $x\in 2a_1B_{(m,k)}$ and $y\in a_2B_{(m,k)}$,
we have $\rz(x)\sim \rz(x_{(m,k)})\sim \rz(y)$, and
there is  a constant $a_3>1$ independent of  $B_{(m,k)}$
such that
\begin{align}
I_2&\lesssim
\int_{a_2B_{(m,k)}}
\l(\int_{|x-y|<a_3\rz(y)}
\frac{\bar{\wz}(x)}{(\rz(y))^n}\,dx\r)
|f\chi_{a_2B_{(m,k)}}(y)|\,dy \nonumber \\
&\lesssim  \int_{a_2B_{(m,k)}} |f(y)| \bar{\wz}(y)\,dy \nonumber\\
&\lesssim \l\|f\r\|_{L^1_{\wz}(a_2B_{(m,k)})}.
\end{align}
We now estimate $I_3$.
For $x\in 2a_1B_{(m,k)}$ and $y\in a_2B_{(m,k)}$,
we have $\rz(x)\sim \rz(x_{(m,k)})\sim \rz(y)$, and
there exists a constant $a_4>1$ independent of  $B_{(m,k)}$
such that $\rz(x)\le a_4\rz(y)$.
On the other hand, if $V\in RH_n(\rn)$,
then there exists $\varepsilon>0$, which
depends only on $n$ and the constant $C$ in (2.1),
such that $V\in RH_{n+\varepsilon}(\rn)$ (see \cite{Ge}).
Thus, according to Lemma 3 of \cite{BHS0}, we get
\beqs
\begin{aligned}
I_3
&\lesssim \int_{2a_1B_{(m,k)}}\int_{\rn}
|R_j(x,y)-\widetilde{R}_j(x,y)|\l|\eta\l(\frac{|x-y|}{\rz(x)}\r)\r|
|f\chi_{a_2B_{(m,k)}}(y)|\,dy\,\bar{\wz}(x)\,dx\\
&\lesssim \int_{2a_1B_{(m,k)}}\int_{|x-y|<\rz(x)}
|R_j(x,y)-\widetilde{R}_j(x,y)|
|f\chi_{a_2B_{(m,k)}}(y)|\,dy\,\bar{\wz}(x)\,dx\\
&\lesssim \int_{2a_1B_{(m,k)}}\int_{|x-y|<\rz(x)}
\frac{|f\chi_{a_2B_{(m,k)}}(y)|}{|x-y|^n}
\l(\frac{|x-y|}{\rz(x)}\r)^{2-n/(n+\varepsilon)}
\,dy\,\bar{\wz}(x)\,dx\\
&\lesssim \int_{a_2B_{(m,k)}}
\l(\int_{|x-y|<a_4\rz(y)}
\frac{1}{|x-y|^n}
\l(\frac{|x-y|}{\rz(y)}\r)^{2-n/(n+\varepsilon)}
\bar{\wz}(x)\,dx\r)
|f(y)|\,dy.\\
\end{aligned}
\eeqs
For the inner integral, there exists  a constant $a_5$ independent of  $B_{(m,k)}$
such that $a_5>a_2$ and $B(y,\max\{a_3,\,2a_4\}\rz(y))\subset a_5B_{(m,k)}$
for all $y\in a_2B_{(m,k)}$, then by the properties of $A_1$ weights, we have
\beqs
\begin{aligned}
&\int_{|x-y|<a_4\rz(y)}
\frac{1}{|x-y|^n}
\l(\frac{|x-y|}{\rz(y)}\r)^{2-n/(n+\varepsilon)}
\bar{\wz}(x)\,dx\\
&\lesssim \sum_{k=k_0}^{\fz}\int_{|x-y|\sim 2^{-k}\rz(y)}
\frac{1}{|x-y|^n}
\l(\frac{|x-y|}{\rz(y)}\r)^{2-n/(n+\varepsilon)}
\bar{\wz}(x)\,dx\\
&\lesssim \sum_{k=k_0}^{\fz}\int_{|x-y|\sim 2^{-k}\rz(y)}
\frac{(2^{-k})^{2-n/(n+\varepsilon)}}{|2^{-k}\rz(y)|^n}
\bar{\wz}(x)\,dx\\
&\lesssim \sum_{k=k_0}^{\fz}
(2^{-k})^{2-n/(n+\varepsilon)}
\int_{|x-y|\le 2^{-k}\rz(y)}
\frac{\bar{\wz}(x)}{|2^{-k}\rz(y)|^n}\,dx\\
&\lesssim  \bar{\wz}(y),
\end{aligned}
\eeqs
where $k_0$ satisfies $2^{-k_0-1}<a_5\le 2^{-k_0}$.
Hence, we obtain
\beq
I_3
\lesssim  \int_{a_2B_{(m,k)}} |f(y)| \bar{\wz}(y)\,dy
\lesssim \l\|f\r\|_{L^1_{\wz}(a_2B_{(m,k)})}.
\eeq
Combining  (3.2)-(3.5), we get
$$\Big\|\sup_{0<2^{-l}<\rz(x)}|\varphi_l*f|\Big\|_{L^1_{\wz}(B_{(m,k)})}
\lesssim \l\|f\r\|_{L^1_{\wz}(a_2B_{(m,k)})}
+\sum_{j=1}^n\|\widehat{R}_j(f)\|_{L^1_{\wz}(2a_1B_{(m,k)})}.
$$
Therefore, by Theorem B and Lemma 2.2, we obtain
\beqs
\begin{aligned}
\|f\|_{h^1_{\rz}(\wz)}
&\lesssim \sum_{(m,k)}
\Big\|\sup_{0<2^{-l}<\rz(x)}|\varphi_l*f|\Big\|_{L^1_{\wz}(B_{(m,k)})}\\
&\lesssim \sum_{(m,k)}
\Big(\l\|f\r\|_{L^1_{\wz}(a_2B_{(m,k)})}
+\sum_{j=1}^n\|\widehat{R}_j(f)\|_{L^1_{\wz}(2a_1B_{(m,k)})}\Big)\\
&\lesssim \|f\|_{L^1_{\wz}(\rn)}
+\sum_{j=1}^n\|\widehat{R}_j(f)\|_{L^1_{\wz}(\rn)}.
\end{aligned}
\eeqs

In order to prove the converse inequality, by Theorem D, it suffices to
show that for all $j$ and any
 $(1,2,s)_{\wz}$-atom or $(1,2)_{\wz}$-single-atom $a$,
$$\|\widehat{R}_j(a)\|_{L^1_{\wz}(\rn)}\lesssim 1.$$
If $a$ is a single atom, by H\"older inequality and
$L^2_{\wz}(\rn)$ boundedness of $\widehat{R}_j$,
we have
$$
\|\widehat{R}_j(a)\|_{L^1_{\wz}(\rn)}
\lesssim \|\widehat{R}_j(a)\|_{L^2_{\wz}(\rn)}\,
\wz(\rn)^{1/2}
\lesssim \|a\|_{L^2_{\wz}(\rn)}\,
\wz(\rn)^{1/2}\lesssim 1.
$$
Next we assume that $a$ is a $(1,2,s)_{\wz}$-atom and
$\supp \, a\subset Q(x_0,r)$ with $r\le L_{1}\rz(x_0)$,
then $\rz(x)\le C_0^22^{k_0}(1+L_1)\rz(x_0)$
for any $x$  satisfying  $|x-y|<\rz(x)$ and $y\in Q(x_0,r)$.

We first consider the atom $a$ with $L_{2}\rz(x_0)\le r\le L_{1}\rz(x_0)$.
Taking $b_1=1+C_0^22^{k_0}(1+L_1)/L_2$,
then $\supp(\widehat{R}_j(a))\subset Q(x_0,b_1r) $,
thus, by H\"older inequality, Lemma 2.4 and
$L^2_{\wz}(\rn)$ boundedness of $\widehat{R}_j$,
we get
\beqs
\begin{aligned}
\|\widehat{R}_j(a)\|_{L^1_{\wz}(\rn)}
&=\|\widehat{R}_j(a)\|_{L^1_{\wz}(Q(x_0,b_1r))}
\lesssim \|a\|_{L^2_{\wz}(\rn)}\,
\wz(Q(x_0,b_1r))^{1/2}\\
&\lesssim \|a\|_{L^2_{\wz}(\rn)}\,
\wz(Q(x_0,r))^{1/2}
\lesssim 1.
\end{aligned}
\eeqs

If  $r< L_{2}\rz(x_0)$,
then $\supp(\widehat{R}_j(a))\subset Q(x_0,b_2\rz(x_0))$,
where $b_2=C_0^22^{k_0}(1+L_1)+L_2$.
Applying the H\"older inequality, Lemma 2.4 and
$L^2_{\wz}(\rn)$ boundedness of $\widehat{R}_j$
give us that
$$\|\widehat{R}_j(a)\|_{L^1_{\wz}(Q(x_0,2r))}
\lesssim \|a\|_{L^2_{\wz}(\rn)}\,
\wz(Q(x_0,r))^{1/2}
\lesssim 1.$$
Moreover, for $x\in Q(x_0,b_2\rz(x_0))\backslash Q(x_0,2r)$,
there is a constant $b_3>1$ such that
 $b_3^{-1}\rz(x_0)\le \rz(x)\le b_3\rz(x_0)$,
 and by the vanishing moment of $a$, we get
\beqs
\begin{aligned}
\l|\widehat{R}_j(a)(x)\r|
&\lesssim \l|\int_{\rn}\frac{x_j-y_j}{|x-y|^{n+1}}
\l[\eta \l(\frac{|x-y|}{\rz(x)}\r)-1\r]a(y)\,dy\r|\\
&\quad + \l|\int_{\rn}
\l[\frac{x_j-y_j}{|x-y|^{n+1}}-\frac{x_j-(x_0)_j}{|x-x_0|^{n+1}}\r]a(y)\,dy\r|\\
&\lesssim \int_{|x-y|\ge \rz(x_0)/2b_3}
\frac{|a(y)|}{|x-y|^n}\,dy
+\int_{\rn}\frac{|x_0-y|}{|x-x_0|^{n+1}}|a(y)|\,dy\\
&\lesssim \l[\frac{1}{(\rz(x_0))^n}+\frac{r}{|x-x_0|^{n+1}}\r]
\|a\|_{L^1(\rn)},
\end{aligned}
\eeqs
where in the penultimate inequality, we use the fact that if
$|x-y|<\rz(x_0)/2b_2$, then $\eta (|x-y|/\rz(x))=1$.
This implies that
\beqs
\begin{aligned}
\int_{Q(x_0,2r)^{\complement}}
\l|\widehat{R}_j(a)(x)\r|\wz(x)\,dx
&=\int_{Q(x_0,b_2\rz(x_0))\backslash Q(x_0,2r)}
\l|\widehat{R}_j(a)(x)\r|\wz(x)\,dx\\
&\lesssim
\|a\|_{L^1(\rn)}
\int_{Q(x_0,b_2\rz(x_0))\backslash Q(x_0,2r)}
\frac{\wz(x)}{(\rz(x_0))^n}\,dx\\
&\quad
+\|a\|_{L^1(\rn)}
\int_{Q(x_0,b_2\rz(x_0))\backslash Q(x_0,2r)}
\frac{r\wz(x)}{|x-x_0|^{n+1}}\,dx\\
&\equiv I+II.
\end{aligned}
\eeqs
For $I$, since $\supp\, a\subset Q_0\equiv Q(x_0,r)\subset Q(x_0,\rz(x_0))\equiv Q$,
by Lemma 2.4 and 2.5, we have
$$I\lesssim \|a\|_{L^1(\rn)}
\frac{\wz(Q(x_0,b_2\rz(x_0)))}{(\rz(x_0))^n}
\lesssim
\frac{|Q|}{\wz(Q)}
\frac{\wz(Q(x_0,b_2\rz(x_0)))}{(\rz(x_0))^n}
\lesssim 1.$$
For $II$, by Lemma 2.4 and 2.5, we get
\beqs
\begin{aligned}
II&\lesssim \|a\|_{L^1(\rn)}
\int_{2r\le|x-x_0|\le b_2\rz(x_0)}
\frac{r\wz(x)}{|x-x_0|^{n+1}}\,dx\\
&\lesssim \|a\|_{L^1(\rn)}
\sum_{j=2}^{j_0}
\int_{|x-x_0|\sim 2^jr}
\frac{r\wz(x)}{|x-x_0|^{n+1}}\,dx\\
&\lesssim \frac{|Q_0|}{\wz(Q_0)}
\sum_{j=2}^{j_0} 2^{-k}
\frac{\wz(2^jQ_0)}{|2^kQ_0|}
\lesssim 1,
\end{aligned}
\eeqs
where $j_0$ is an integer such that $2^{j_0-1}r<b_2\rz(x_0)\le 2^{j_0}r$.
Thus, for $r< L_{2}\rz(x_0)$ we obtain
$$\|\widehat{R}_j(a)\|_{L^1_{\wz}(\rn)}\lesssim 1.$$

 Finally, we prove
 $\|f\|_{L^1_{\wz}(\rn)}
+\sum_{j=1}^n\|{R}_j(f)\|_{L^1_{\wz}(\rn)}\sim
\|f\|_{L^1_{\wz}(\rn)}
+\sum_{j=1}^n\|\widehat{R}_j(f)\|_{L^1_{\wz}(\rn)}.$ Indeed, we only need to show that
\beqs
\|\widehat{R}_j(f)-{R}_j(f)\|_{L^1_\wz(\rn)}\lesssim\|f\|_{L^1_\wz(\rn)}.
\eeqs
In fact, note that $|x-y|\le \rho(x)$, then there exists a constant $\bz>0$ such that $\rho(y)/\bz\le \rho(x)\le \bz\rho(y)$ by Lemma 2.1. Set $\wz\in A_1^{\rho,\tz}(\rn)$ with $\tz>0$. Hence, by the properties of $A_1^{\rho,\tz}(\rn)$, we get
\beqs
\begin{aligned}
\|\widehat{R}_j(f)-{R}_j(f)\|_{L^1_\wz(\rn)}&\le \int_{\rn}\int_{\rn}|R_j(x,y)-\widetilde{R}_j(x,y)|\eta\l(\frac {|x-y|}{\rho(x)}\r)|f(y)|dy\wz(x)dx\\
&\le  \int_{\rn}\int_{|x-y|<\rho(y)/\bz}|R_j(x,y)-\widetilde{R}_j(x,y)|\wz(x)dx|f(y)|dy\\
&\qquad+\int_{\rn}\int_{\rho(y)/\bz\le |x-y|<\bz\rho(y)}|\widetilde{R}_j(x,y)|\wz(x)dx|f(y)|dy\\
&\qquad+\int_{\rn}\int_{\rho(y)/\bz\le |x-y|}|{R}_j(x,y)|\wz(x)dx|f(y)|dy\\
&\lesssim \|f\|_{L^1_\wz(\rn)}+\int_{\rn} \frac{\wz(Q(y,\rz(y)))}{(\rz(y))^n}|f(y)|dy\\
&\qquad+\int_{\rn}\dsum_{k=1}^\fz 2^{-kN}|2^k\rho(y)|^{-n}\int_{ |x-y|\le 2^k\rho(y)/\bz}\wz(x)dx|f(y)|dy\\
&\lesssim \|f\|_{L^1_\wz(\rn)}+\int_{\rn}\dsum_{k=1}^\fz 2^{-k(N-(k_0+1)\tz)}\frac{\wz(Q(y,\rz(y)))}{(\rz(y))^n}|f(y)|dy\\
&\lesssim\|f\|_{L^1_\wz(\rn)},
\end{aligned}
\eeqs
if taking $N>(k_0+1)\tz+1$.

Therefore, the proof is complete.
\qed

\bigskip

\noi{\bfseries{Acknowledgements}}\quad
The research is supported by National Natural Science Foundation of China under Grant
\#11271024.

\bigskip
\bigskip

\indent Beijing International Studies University,\\
\indent Beijing, 100024, \\
\indent People’s Republic of China

\textsl{E-mail address: zhuhua@pku.edu.cn}


\begin{thebibliography}{99}



\bibitem{BHS0}  \label{BHS0} B. Bongioanni, E. Harboure, O. Salinas,
{\it Riesz transforms related to Schr\"odinger operators acting
on $BMO$ type spaces,}
J. Math. Anal. Appl. $\mathbf{357}$ (2009), 115-131.


\bibitem{BHS1}  \label{BHS1} B. Bongioanni, E. Harboure, O. Salinas,
{\it Classes of weights related to Schr\"odinger operators,}
J. Math. Anal. Appl. $\mathbf{373}$ (2011), 563-579.

\bibitem{BHS2}  \label{BHS2} B. Bongioanni, E. Harboure, O. Salinas,
{\it Commutators of Riesz transforms related to Schr\"odinger operators,}
J. Fourier Anal. Appl. $\mathbf{17}$ (2011), 115-134.


\bibitem{BLYZ}  \label{BLYZ} M. Bownik, B. Li, D. Yang and Y. Zhou,
{\it Weighted anisotropic Hardy spaces and their
applications in boundedness of sublinear operators,}
Indiana Univ. Math. J. $\mathbf{ 57}$ (2008), 3065-3100.


\bibitem{Du}  \label{Du}  J. Duoandikoetxea,
{\it Fourier Analysis,}
Grad. Stud. Math., vol. 29, Amer. Math. Soc., Providence, RI, 2000.


\bibitem{DG} \label{DG} J. Dziuba\'{n}ski, G. Garrig\'{o}s,
T. Mart\'{i}nez, J. L. Torrea and J. Zienkiewicz,
{\it BMO spaces related to Schr\"{o}dinger operators with
potentials satisfying a reverse H\"{o}der inequality,}
Math. Z. $\mathbf{249}$ (2005), 329-356.


\bibitem{DZ1}  \label{DZ1} J. Dziuba\'{n}ski and J. Zienkiewicz,
{\it Hardy space $H^{1}$ associated to Schr\"{o}dinger operator with
potential satisfying reverse H\"{o}lder inequality,}
Rev. Mat. Iberoam. $\mathbf{15}$ (2) (1999), 279-296.

\bibitem{DZ2}  \label{DZ2} J. Dziuba\'{n}ski and J. Zienkiewicz,
{\it $H^{p}$ spaces for Schr\"{o}dinger operators,}
Fourier Analysis and Related Topics, Banach Center Publications. $\mathbf{56}$ (2002), 45-53.

\bibitem{DZ3}  \label{DZ3} J. Dziuba\'{n}ski and J. Zienkiewicz,
{\it $H^{p}$ spaces associated with Schr\"{o}dinger operators
with potentials from reverse H\"{o}lder classes,}
Colloq. Math. $\mathbf{98}$ (2003), 5-38.


\bibitem{FS} \label{FS} C. Fefferman and E. Stein,
{\it $H^{p}$ spaces of several variables,}
Acta. Math. $\mathbf{129}$ (1972), no. 3-4, 137-193.

\bibitem{Ga} \label{Ga} J. Garc\'ia-Cuerva,
{\it Weighted $H^{p}$ spaces,}
Dissertationes Math. (Rozprawy Mat.) $\mathbf{162}$ (1979), 63pp.


\bibitem{GR} \label{GR} J. Garc\'ia-Cuerva and J. Rubio de Francia,
{\it Weighted norm inequalities and related topics,}
North-Holland, Amsterdam.  1985.

\bibitem{Ge} \label{Ge} F. Gehring,
{\it The $L^{p}$-integrability of the partial
derivatives of a quasi-conformal mapping,}
Acta Math. $\mathbf{130}$ (1973), 265-277.


\bibitem{Go} \label{Go} D. Goldberg,
{\it A local version of real Hardy spaces,}
Duke Math. J. $\mathbf{46}$ (1979), 27-42.


\bibitem{Gr} \label{Gr} L. Grafakos,
{\it Classical and Modern Fourier Analysis,}
Pearson Edu. Inc., New Jersey, 2004.


\bibitem{LTZ} \label{LTZ} H. Liu, L. Tang and H. Zhu,
{\it Weighted Hardy spaces and $BMO$ spaces
associated with Schr\"odinger operators,}
Math. Nachr.  $\mathbf{285}$ (2012), 2173-2207.

\bibitem{MPR} \label{MPR} G. Mauceri, M. Picardello and F. Ricci,
{\it A Hardy space associated with twisted convolution,}
Adv. Math.  $\mathbf{39}$ (1981), 270-288.






\bibitem{Mu} \label{Mu} B. Muckenhoupt,
{\it Weighted norm inequalities for the Hardy maximal functions,}
Trans. Amer. Math. Soc.   $\mathbf{165}$ (1972), 207-226.


\bibitem{Ry} \label{Ry} V. S. Rychkov,
{\it Littlewood-Paley theory and function spaces with $A_p^{loc}$ weights,}
Math. Nachr.  $\mathbf{224}$ (2001), 145-180.



\bibitem{Sh} \label{Sh} Z. Shen,
{\it $L^{p}$ estimates for Schr\"{o}dinger operators with certain potentials,}
Ann. Inst. Fourier (Grenoble) $\mathbf{45}$ (2) (1995), 513-546.



\bibitem{St} \label{St} E. M. Stein,
{\it Harmonic Analysis: Real Variable Methods,
Orthogonality, and Oscillatory Integrals,}
Princeton University Press, Princeton, N.J., 1993.


\bibitem{ST} \label{ST} J. Stromberg and A. Torchinsky,
{\it Weighted Hardy spaces,}
Springer-Verlag, Berlin, 1989.


\bibitem{Ta1} \label{Ta1} L. Tang,
{\it Weighted local Hardy spaces and their applications,}
Illinois J. Math. $\mathbf{56}$(2012), 453-495.


\bibitem{Ta2} \label{Ta2} L. Tang,
{\it Weighted norm inequalities for pseudo-differential
operators with smooth symbols and their commutators,}
J. Funct. Anal. $\mathbf{261}$(2012), 1602-1629.

\bibitem{Ta3} \label{Ta3} L. Tang,
{\it Weighted norm inequalities for Schr\"odinger type operators,}
to be appeared in Forum Math. DOI 10.1515/forum-2013-0070(arXiv:1109.0099 v1).

\bibitem{Ta4} \label{Ta4} L. Tang,
{\it  Extrapolation from $A_\fz^{\rho,\fz}$, vector-valued
inequalities and applications in the  Schr\"odinger settings,}
to be appeared in Ark. Mat. (arXiv:1109.0101v1).





\bibitem{YYZ} \label{YYZ} D. Yang, D. Yang and Y. Zhou,
{\it Endpoint properties of localized Riesz transforms and fractional
integrals associated to Schr\"odinger operators,}
Potential Anal.  $\mathbf{30}$ (2009), 271-300.

\bibitem{YZ} \label{YZ} D. Yang and Y. Zhou,
{\it Boundedness of sublinear operators in Hardy spaces on RD-spaces via atoms,}
J. Math. Anal. Appl. $\mathbf{339}$ (2008), 622-635.

\bibitem{YZ1} \label{YZ1} D. Yang and Y. Zhou,
{\it Localized Hardy spaces $H^1$ related to admissible functions on RD-spaces and applications to Schr\"odinger operators,}
Trans. Amer. Math. Soc. $\mathbf{363}$ (2011), 1197-1239.

\bibitem{Zh} \label{Zh} J. Zhong,
{\it Harmonic analysis for some Schr\"odinger type operators,}
Ph. D. Thesis. Princeton University, 1993.

\bibitem{ZL} \label{ZL} H. Zhu and H. Liu,
{\it Weighted estimates for bilinear operators,}
J. Funct. Spaces Appl., vol. 2014, Article ID 797956, 10 pages, 2014.


\bibitem{ZT} \label{ZL} H. Zhu and L. Tang,
{\it Weighted local Hardy spaces associated to Schr\"{o}dinger operators,}
arXiv:1403.7641.





\end{thebibliography}
\end{document}